\documentclass[12pt]{amsart}
\usepackage{amsxtra,latexsym,amssymb}
\usepackage{epsfig,rotate,amsthm}
\usepackage{hyperref}

\def\qed{{\hfill $\Box$}}

\def\Z{{\mathbb Z}}
\def\C{{\mathbb C}}

\def\U{{U_{r,s}^{+}({\mathfrak sl_{3}})}}
\def\V{{\check U}^{\geq 0}_{r,s}({\mathfrak sl_{3}})}
\theoremstyle{theorem}
\newtheorem{thm}{Theorem}[section]
\newtheorem{cor}{Corollary}[section]
\newtheorem{prop}{Proposition}[section]

\newtheorem{lem}{Lemma}[section]
\theoremstyle{definition}
\newtheorem{defn}{Definition}[section]
\theoremstyle{remark}

\begin{document}
\title[(Hopf) Algebra Automorphisms]{(Hopf) Algebra Automorphisms of
the Hopf algebra ${\check U}^{\geq 0}_{r,s}({\mathfrak sl_{3}})$}
\author[X. Tang]{Xin Tang}
\address{Department of Mathematics \& Computer Science\\
Fayetteville State University\\
1200 Murchison Road, Fayetteville, NC 28301}
\email{xtang@uncfsu.edu} 
\thanks{This research project is partially supported by Funds from the HBCU
Master's Degree STEM Program at Fayetteville State University.}
\date{\today}
\begin{abstract}
In this paper, we completely determine the group of algebra
automorphisms for the two-parameter Hopf algebra ${\check U}_{r,s}^{\geq
0}({\mathfrak sl_{3}})$. As a result, the group of Hopf 
algebra automorphisms is determined for $\V$ as well. We further 
characterize all the derivations of the subalgebra $U^{+}_{r,s}({\mathfrak
sl_{3})}$, and calculate its first degree Hochschild cohomology group.
\end{abstract}
\maketitle 
\section*{Introduction}
Motivated by the study of down-up algebras \cite{BT}, a two-parameter quantized enveloping 
algebra (quantum group) $U_{r,s}({\mathfrak sl_{n}})$ has recently been investigated by 
Benkart and Witherspoon in the references \cite{BW1, BW2}. The two-parameter quantized enveloping algebra $U_{r,s}({\mathfrak sl_{n}})$ remains as a close analogue of the one-parameter quantized enveloping algebra $U_{q}({\mathfrak sl_{n}})$ associated to the finite dimensional complex simple special linear Lie algebra ${\mathfrak sl_{n}}$. As a matter of fact, 
the two-parameter quantized enveloping algebra $U_{r,s}({\mathfrak sl_{n}})$ shares many 
similar properties with its one-parameter analogue. 

For instance, the two-parameter quantized enveloping algebra $U_{r,s}({\mathfrak sl_{n}})$ is also a Hopf algebras and it admits a triangular decomposition. Besides having a similar representation theory to that of the one-parameter quantum group $U_{q}({\mathfrak sl_{n}})$, the Hopf algebra $U_{r,s}({\mathfrak sl_{n}})$ can be realized as the Drinfeld double of its Hopf sub-algebras $U_{r,s}^{\geq 0}({\mathfrak sl_{n}})$ and $U^{\leq 0}_{r,s}({\mathfrak sl_{n}})$. However, the algebra $U_{r,s}({\mathfrak sl_{n}})$ does have some different ring-theoretic features in that $U_{r,s}({\mathfrak sl_{n}})$ is more rigid and possesses less symmetries. In addition, the center of the algebra $U_{r,s}({\mathfrak sl_{n}})$ shows a different picture as well \cite{BKL}.

To better understand the structure and properties of the algebra $U_{r,s}({\mathfrak sl_{n})}$, one naturally starts with the investigations of its subalgebra $U_{r,s}^{+}({\mathfrak sl_{n}})$ and Hopf subalgebra $U^{\geq 0}_{r,s}({\mathfrak sl_{n}})$. In this paper, we will study the algebra $U_{r,s}^{+}({\mathfrak sl_{3}})$ in terms of its derivations, and its augmented Hopf algebra $\V$ in terms of its algebra automorphisms and Hopf algebra automorphisms. We shall determine all the algebra automorphisms and Hopf algebra automorphisms of $\V$. In addition, we will characterize all the derivations of $U_{r,s}^{+}({\mathfrak sl_{3}})$. As a result, we will 
compute the first Hochschild cohomology group of the algebra $U_{r,s}^{+}({\mathfrak sl_{3}})$. 

Now let me briefly mention the methods which we shall follow. In order to determine the algebra and Hopf algebra automorphisms of $\V$, we shall closely follow the approach used in \cite{F}. In order to characterize the derivations of $\U$, we shall embed the algebra $\U$ into a quantum torus, whose derivations had been explicitly described in \cite{OP}. We would like to point out that this embedding will allow us to extend the derivations of $\U$ to the derivations of the associated quantum torus. Therefore, via this embedding, we shall be able to pull the information on derivations back to the algebra $U_{r,s}^{+}({\mathfrak sl_{3}})$. Based on a result on the derivations of quantum torus established in \cite{OP}, we will be able to determine all the derivations of the algebra $U^{+}_{r,s}({\mathfrak sl_{3}})$ modulo its inner derivations. As an immediate application, we show that the first Hochschild cohomology group $HH^{1}(\U)$ of $\U$ is a $2-$dimensional vector space over the base field $\C$.

The paper is organized as follows. In Section 1, we recall some basic definitions and properties on 
the two-parameter quantized enveloping algebras $U^{+}_{r,s}({\mathfrak sl_{3}})$ and its 
augmented Hopf algebra $\V$. In Section 2, we determine the algebra automorphism group 
and Hopf algebra automorphism group of $\V$. In Section 3, we characterize the derivations of 
$U_{r,s}^{+}({\mathfrak sl_{3}})$, and compute the first Hochschild cohomology group 
$HH^{1}(U_{r,s}^{+}({\mathfrak sl_{3}}))$.

\section{Definitions and Basic Properties of $U_{r,s}^{+}({\mathfrak sl_{3}})$ and ${\check U}^{\geq 0}_{r,s}({\mathfrak sl_{3}})$}

\subsection{Definition and basic properties of $U_{r,s}^{+}({\mathfrak sl_{3}})$}
Motivated by the study of down-up algebras, a two-parameter quantized enveloping algebra 
$U_{r,s}({\mathfrak sl_{n}})$ associated to the finite dimensional simple complex Lie algebra 
${\mathfrak sl_{n}}$ has been recently studied by Benkart and Witherspoon in \cite{BW1} 
and the references therein. For the purpose of this paper, we shall only recall the definitions 
of the algebras $U^{+}_{r,s}({\mathfrak sl_{n}})$ and $U_{r,s}^{\geq 0}({\mathfrak sl_{n}})$, 
which are indeed subalgebras of $U_{r,s}({\mathfrak sl_{n}})$. One easily sees that the algebra 
$U^{+}_{r,s}({\mathfrak sl_{n}})$ can be regarded as a two-parameter quantized enveloping 
algebra of a maximal nilpotent Lie subalgebra of the Lie algebra ${\mathfrak sl_{n}}$. 

Let $C=(a_{ij})$ denote the Cartan matrix associated to the Lie algebra $\mathfrak sl_{n}$. Let us define the following notation:
\begin{eqnarray*}
<i,j>= a_{ij} \,\text{for}\, i<j;\\
<i,i>=1\,\text{for}\, i=1, \cdots, n-1;\\
<i,j>=0 \, \text{for}\, i>j.
\end{eqnarray*}

Suppose that $r,s\in \C$ such that $r^{m}s^{n}=1$ implies $m=n=0$. We recall the following definition:

\begin{defn}
The two-parameter quantized enveloping algebra $U_{r,s}^{\geq 0}({\mathfrak sl_{n}})$ is defined to be the $\C-$algebra 
generated by $E_{i}, W_{i}$ subject to the following relations:
\begin{eqnarray*}
W_{i}^{\pm 1}W_{j}^{\pm 1}=W_{j}^{\pm 1}W_{i}^{\pm 1};\\
W_{i}^{\pm 1}W_{i}^{\mp 1}=1;\\
W_{i}E_{j}=r^{<j,i>}s^{-<i,j>}E_{j}W_{i};\\
E_{i}^{2}E_{i+1}-(r+s)E_{i}E_{i+1}E_{i}+rsE_{i+1}E_{i}^{2}=0;\\
E_{i+1}^{2}E_{i}-(r^{-1}+s^{-1})E_{i+1}E_{i}E_{i+1}+r^{-1}s^{-1}E_{i}E_{i+1}^{2}=0.\\
\end{eqnarray*}
And the algebra $U_{r,s}^{+}({\mathfrak sl_{n}})$ is defined to be the subalgebra of $U_{r,s}^{\geq 0}({\mathfrak sl_{n}})$ generated by $E_{i}$.
\end{defn}

From \cite{BW1}, we know that the two-parameter quantized enveloping algebra $U_{r,s}^{\geq 0}({\mathfrak sl_{n}})$ has a Hopf algebra structure, which is defined by the following coproduct, counit and antipode:
\begin{eqnarray*}
\Delta(W_{i}^{\pm 1})=W_{i}^{\pm 1}\otimes W_{i}^{\pm 1};\\
\Delta(E_{i})=E_{i}\otimes 1+ W_{i}\otimes E_{i}; \\
\epsilon ( W_{i}^{\pm 1})=1;\\
\epsilon (E_{i})=0;\\
S(W_{i}^{\pm 1})=W_{i}^{\mp 1};\\
S(E_{i})=-W_{i}^{-1}E_{i}.
\end{eqnarray*}

When $n=3$, one has the corresponding two-parameter quantized enveloping algebra $U_{r,s}^{+}({\mathfrak sl_{3}})$, which will be of one of the major objects in this paper. 
In particular, we recall the following definition:
\begin{defn}
The algebra $U^{+}_{r,s}({\mathfrak sl_{3}})$ is defined to be the $\C-$algebra generated by the generators $E_{1}, E_{2}$ subject to 
the following relations:
\begin{eqnarray*}
E_{1}^{2}E_{2}-(r+s)E_{1}E_{2}E_{1}+rsE_{2}E_{1}^{2}=0;\\
E_{1}E_{2}^{2}-(r+s)E_{1}E_{2}E_{1}+rsE_{2}^{2}E_{1}=0.
\end{eqnarray*}
\end{defn}

Naturally, one may think of the algebra $U_{r,s}^{+}({\mathfrak sl_{3}})$ as a two-parameter quantum Heisenberg algebra. Indeed, the algebra $U_{r,s}^{+}({\mathfrak sl_{3}})$ shares many similar properties as the algebra $U_{q}^{+}({\mathfrak sl_{3}})$, which has been traditionally called the quantum Heisenberg algebra. 

In addition, we recall the definition of the following Hopf subalgebra of $U_{r,s}({\mathfrak sl_{3}})$:
\begin{defn}
The Hopf algebra $U_{r,s}^{\geq 0}({\mathfrak sl_{3}})$ is defined to be the $\C-$algebra generated by $E_{1}, E_{2}, W_{1}, W_{2}$ subject to the following relations:
\begin{eqnarray*}
W_{1}W_{1}^{-1}=1=W_{2}W_{2}^{-1};\\
W_{1}W_{2}=W_{2}W_{1};\\
W_{1}E_{1}=rs^{-1}E_{1}W_{1};\\
W_{1}E_{2}=sE_{2}W_{1};\\
W_{2}E_{1}=r^{-1}E_{1}W_{2};\\
W_{2}E_{2}=rs^{-1}E_{2}W_{2};\\
E_{1}^{2}E_{2}-(r+s)E_{1}E_{2}E_{1}+rsE_{2}E_{1}^{2}=0;\\
E_{1}E_{2}^{2}-(r+s)E_{2}E_{1}E_{2}+rsE_{2}^{2}E_{1}=0.
\end{eqnarray*}
\end{defn}

However, we should mention that we will not study the Hopf algebra $U_{r,s}^{\geq 0}({\mathfrak sl_{3}})$ in this paper. Instead, we will study its augmented version $\V$, which we shall define in the next subsection. 

Before we introduce the Hopf algebra $\V$, let us mention some basic properties of the algebra $\U$ in the rest of this subsection. It is easy to see that the two-parameter quantized enveloping algebra $U_{r,s}^{+}({\mathfrak sl_{n}})$ can also be presented as an iterated skew polynomial ring and a PBW-basis can be constructed for $U_{r,s}({\mathfrak sl_{n}})$ as well. For conveniences, we shall only recall the skew polynomial presentation for the algebra $\U$ here. 
For the general construction, we refer the reader to references \cite{BKL,T}. 

First of all, Let us fix some notation by setting the following new variables: 
\[
E_{1}=E_{1}, \quad E_{2}=E_{2},\quad E_{3}=E_{1}E_{2}-sE_{2}E_{1}.
\]

Then it is easy to see that we have the following relations between these 
new variables:
\[
E_{1}E_{3}=rE_{3}E_{1},\quad E_{2}E_{3}=r^{-1}E_{3}E_{2}.
\]

Now let us further define some algebra automorphisms $\tau_{2},\tau_{3}$ and some derivations $\delta_{2},\delta_{3}$ as follows: 
\begin{eqnarray*}
\tau_{2}(E_{1})=r^{-1}E_{1};\\
\delta_{2}(E_{1})=0;\\
\tau_{3}(E_{1})=s^{-1}E_{1};\\
\tau_{3}(E_{3})=r^{-1}E_{3};\\
\delta_{3}(E_{1})=-s^{-1}E_{3};\\
\delta_{3}(E_{3})=0.
\end{eqnarray*}
\qed

Then it is easy to see that we have the following result
\begin{thm}
The algebra $\U$ can be presented as an iterated skew polynomial ring. In
particular, we have the following result
\[
\U \cong
\C[E_{1}][E_{3},\tau_{2},\delta_{2}][E_{2},\tau_{3},\delta_{3}].
\] 
\end{thm}
\qed

Based on the previous theorem, we have an obvious corollary as follows:
\begin{cor}
The set $\{E_{1}^{i}E_{3}^{j}E_{2}^{k}|i, j, k\geq 0\}$ forms 
a PBW-basis of the algebra $\U$. In particular, $\U$ has a $GK-$dimension 
of $3$.
\end{cor}
\qed

\subsection{The Augmented Hopf algebra ${\check U}^{\geq 0}_{r,s}({\mathfrak sl_{3}})$}
In this subsection, we shall introduce an augmented Hopf algebra $\V$, which contains the algebra 
$U^{+}_{r,s}({\mathfrak sl_{3}})$ as a subalgebra and enlarges the Hopf algebra $U^{\geq 0}_{r,s}({\mathfrak (sl_{3}})$.

First of all, we need to define the following new variables:
\[
K_{1}=W_{1}^{2/3}W_{2}^{1/3},\quad K_{2}=W_{1}^{1/3}W_{2}^{2/3}\\
\]

Then we have the following definition of $\V$.
\begin{defn}
The algebra ${\check U}^{\geq 0}_{r,s}({\mathfrak sl_{3}})$ is a $\C$-algebra generated by 
$E_{1}, E_{2}, K_{1}^{\pm 1}, K_{2}^{\pm 1}$ subject to the following relations:
\begin{eqnarray*}
K_{1}K_{1}^{-1}=1=K_{2}K_{2}^{-1};\\
K_{1}K_{2}=K_{2}K_{1};\\
K_{1}E_{1}=r^{1/3}s^{-2/3}E_{1}K_{1};\\
K_{1}E_{2}=r^{1/3}s^{1/3} E_{2}K_{1};\\
K_{2}E_{1}=r^{-1/3}s^{-1/3}E_{1}K_{2};\\
K_{2}E_{2}=r^{2/3}s^{-1/3}E_{2}K_{2};\\
E_{1}^{2}E_{2}-(r+s)E_{1}E_{2}E_{1}+rsE_{2}E_{1}^{2}=0;\\
E_{1}E_{2}^{2}-(r+s)E_{2}E_{1}E_{2}+rsE_{2}^{2}E_{1}=0.
\end{eqnarray*}
\end{defn}
\qed

In order to introduce a Hopf algebra structure on $\V$, let us further define some operators as follows:
\begin{eqnarray*}
\Delta(E_{1})=E_{1}\otimes 1+K_{1}^{2}K_{2}^{-1}\otimes E_{1};\\
\Delta(E_{2})=E_{2}\otimes 1+K_{1}^{-1}K_{2}^{2}\otimes E_{2};\\
\Delta(K_{1})=K_{1}\otimes K_{1};\\
\Delta(K_{2})=K_{2}\otimes K_{2};\\
S(E_{1})=-K_{1}^{2}K_{2}^{-1}E_{1};\\
S(E_{2})=-K_{1}^{-1}K_{2}^{2}E_{1};\\
S(K_{1})=K_{1}^{-1};\\
S(K_{2})=K_{2}^{-1};\\
\epsilon(E_{1})=\epsilon(E_{2})=0;\\
\epsilon (K_{1})=\epsilon(K_{2})=1.
\end{eqnarray*}

Then it is straightforward to verify the following result:
\begin{prop}
The algebra ${\check U}^{\geq 0}_{r,s}({\mathfrak sl_{3}})$ is a Hopf
algebra with the coproduct, counit and antipode defined as above.
\end{prop}
\qed

Recall that we have $E_{3}=E_{1}E_{2}-sE_{2}E_{1}$, then it is easy to
see that we have the following result 
\begin{thm}
The algebra ${\check U}^{\geq 0}_{r,s}(sl_{3})$ has a $\C-$basis 
\[
\{ K_{1}^{m}K_{2}^{n}E_{1}^{i}E_{2}^{j}E_{3}^{k}\mid m, n\in \Z,\, i,j,k\in \Z_{\geq 0}\}.
\]
\end{thm}
\qed

In particular, one can see that all the invertible elements of $\V$ are of the form $\lambda K_{1}^{m}K_{2}^{n}$ for some $\lambda \in \C^{\ast}$ and $m, n \in \Z$.

\section{Algebra and Hopf algebra automorphisms of ${\check U}^{\geq 0}_{r,s}(\mathfrak sl_{3})$}

In this section, we will first determine the algebra automorphism group of the algebra ${\check U}_{r,s}^{\geq 0}({\mathfrak sl_{3}})$. As a result, we
are able to determine the Hopf algebra automorphism group of ${\check U}_{r,s}^{\geq 0}({\mathfrak sl_{3}})$ as well. We will closely
follow the approach used in \cite{F}.

\subsection{The algebra automorphism group of $\V$}

Suppose that $\theta \in Aut_{\C}(\V)$ is an algebra automorphism of the algebra $\V$. Note that the elements $K_{1}, K_{2}$ are invertible in $\V$ and $\theta$ is an algebra automorphism, then the elements $\theta(K_{1})$ and $\theta(K_{2})$ are invertible too. Therefore, we have the following 
\[
\theta(K_{1})=a_{1}K_{1}^{x}K_{2}^{y},\quad \theta(K_{2})=a_{2}K_{1}^{z}K_{2}^{w}
\]
for some $a_{1}, a_{2}\in \C^{\ast}$ and $x,y,z,w \in \Z$. 

Let $M_{\theta}=(M_{ij})$ denote the corresponding $2\times 2-$matrix associated to the algebra automorphism $\theta$. Specifically, we will set $M_{11}=x, M_{12}=y, M_{21}=z$ and $M_{22}=w$. Since $\theta$ is an algebra automorphism, we know that the matrix $M_{\theta}$ is 
an invertible matrix with integer coefficients. In particular, we have the following
\[
xw-yz=\pm 1.
\]

Suppose that for $l=1,2$, we have
\[
\theta(E_{l})=\sum_{m_{l}, n_{l}, \beta_{l}^{1}, \beta_{l}^{2}, \beta_{l}^{3}}a_{m_{l}, n_{l}, \beta_{l}^{1}, \beta_{l}^{2}, \beta_{l}^{3}}K_{1}^{m_{l}}K_{2}^{n_{l}}E_{1}^{\beta_{l}^{1}}E_{2}^{\beta_{l}^{2}}E_{3}^{\beta_{l}^{3}}
\]
where $a_{m_{l}, n_{l},\beta_{l}^{1}, \beta_{l}^{2}, \beta_{l}^{3}} \in \C^{\ast}$ and
$m_{l}, n_{l} \in \Z$ and $\beta_{l}^{1}, \beta_{l}^{2}, \beta_{l}^{3} \in \Z_{\geq 0}$. 

Then we have the following
\begin{prop}
Let $\theta \in Aut_{\C}(\V)$ be an algebra automorphism of $\V$, then we have $M_{\theta} \in GL(2,\Z_{\geq 0})$.
\end{prop}
{\bf Proof:} The proof of {\bf Proposition 2.1} in \cite{F} can essentially
be adopted here word by word with the replacement of $s$ by $r^{-1}$. For the reader's convenience, we will present the detailed proof here.

Since $K_{1}E_{1}=r^{1/3}s^{-2/3}E_{1}K_{1}$ and $K_{2}E_{1}=r^{-1/3}s^{-1/3}E_{1}K_{2}$, we have the following
\begin{eqnarray*}
\theta(K_{1})\theta(E_{1})=r^{1/3}s^{-2/3}\theta(E_{1})\theta(K_{1});\\
\theta(K_{2})\theta(E_{1})=r^{-1/3}s^{-1/3}\theta(E_{1})\theta(K_{2}).
\end{eqnarray*}

Via a detailed calculation, we further have the following
\begin{eqnarray*}
(\beta_{1}^{1}+\beta_{1}^{3})x+(\beta_{1}^{2}+\beta_{1}^{3})y=1;\\
(\beta_{1}^{1}+\beta_{1}^{3})z+(\beta_{1}^{2}+\beta_{1}^{3})w=0.
\end{eqnarray*}

Similarly, we also have the following
\begin{eqnarray*}
(\beta_{2}^{1}+\beta_{2}^{3})x+(\beta_{2}^{2}+\beta_{2}^{3})y=0\\
(\beta_{2}^{1}+\beta_{2}^{3})z+(\beta_{2}^{2}+\beta_{2}^{3})w=1.
\end{eqnarray*}

Now let us set a $2\times2-$matrix $B=(b_{ij})$ with the following entries
\begin{eqnarray*}
b_{11}=(\beta_{1}^{1}+\beta_{1}^{3});\\
b_{21}=(\beta_{1}^{2}+\beta_{1}^{3});\\
b_{12}=(\beta_{2}^{1}+\beta_{2}^{3});\\
b_{22}=(\beta_{2}^{2}+\beta_{2}^{3}).
\end{eqnarray*}

Then we have the following system of equations
\begin{eqnarray*}
b_{11}x+b_{21}y=1;\\
b_{12}x+b_{22}y=0;\\
b_{11}z+b_{21}w=0;\\
b_{12}z+b_{22}w=1.
\end{eqnarray*}

This implies that the product $M_{\theta}B$ of matrices $M_{\theta}$ and $B$ is equal to the identity matrix. Therefore, we have $M_{\theta^{-1}}=B$, where $M_{\theta^{-1}}$ is the corresponding matrix associated to the algebra automorphism $\theta^{-1}$. Note that the entries $b_{11}, b_{12}, b_{21}, b_{22}$ are all nonnegative integers, thus we have
$M_{\theta^{-1}}\in GL(2, \Z_{\geq 0})$. 

Applying similar arguments to the algebra automorphism $\theta^{-1}$, we can prove that $M_{\theta}\in GL(2, \Z_{\geq 0})$ as desired. 
\qed

To proceed, we now recall an important lemma ({\bf Lemma 2.2} from \cite{F}), which characterizes 
the matrix $M_{\theta}$:
\begin{lem}
If $M$ is a matrix in $GL(n,\Z_{\geq 0})$ such that it’s inverse matrix $M^{-1}$ is also in $GL(n,\Z_{\geq 0})$, then we 
have $M=(\delta_{i\sigma(j)})_{i,j}$, where $\sigma$ is an element of the symmetric group $\mathbb{S}_{n}$.
\end{lem}
\qed

Based the previous Proposition and Lemma, we immediately have the following result
\begin{cor}
Let $\theta \in Aut_{\C}(\V)$ be an algebra automorphism of $\V$. Then for $l=1,2$, we have 
\[
\theta(K_{l})=a_{l} K_{\sigma(l)}
\]
where $\sigma \in \mathbb{S}_{2}$ and $a_{l} \in \C^{\ast}$.
\end{cor}
\qed

Furthermore, we can further prove the following result:
\begin{prop}
Let $\theta \in Aut_{\C}(\V)$ be an algebra automorphism of $\V$. Then for $l=1,2$, we have 
\[
\theta(E_{l})=b_{l}K_{1}^{m_{l}}K_{2}^{n_{l}}E_{\sigma(l)}
\]
where $b_{l}\in \C^{\ast}$ and $m_{l}, n_{l} \in \Z$.
\end{prop}
{\bf Proof:} Let $\theta \in Aut_{\C}(\V)$ be an algebra automorphism
of $\V$. To prove the proposition, there are two cases to consider:

{\bf Case 1}: Suppose that $\theta(K_{1})=a_{1}K_{1}$ and $\theta(K_{2})=a_{2}K_{2}$, then we need to prove that 
\[
\theta(E_{1})=b_{1}K_{1}^{m_{1}}K_{2}^{n_{1}}E_{1},
\quad
\theta(E_{2})=b_{2}K_{1}^{m_{2}}K_{2}^{n_{2}}E_{2}.
\]

Since $K_{1}E_{1}=r^{1/3}s^{-2/3}E_{1}K_{1}$ and $K_{2}E_{1}=r^{-1/3}s^{-1/3}E_{1}K_{2}$, we have the 
following
\begin{eqnarray*}
\theta(K_{1})\theta(E_{1})=r^{1/3}s^{-2/3}\theta(E_{1})\theta(K_{1});\\
\theta(K_{2})\theta(E_{1})=r^{-1/3}s^{-1/3}\theta(K_{2})\theta(K_{1}).
\end{eqnarray*}

Thus we have the following
\begin{eqnarray*}
&&a_{1}K_{1}(\sum_{m_{1},n_{1},\beta_{1}^{1},\beta_{1}^{2},
\beta_{1}^{3}}b_{m_{1}, n_{1}, \beta_{1}^{1}, \beta_{1}^{2}, \beta_{1}^{3}}
K_{1}^{m_{1}}K_{2}^{n_{1}}E_{1}^{\beta_{1}^{1}}E_{2}^{\beta_{1}^{2}}E_{3}^{\beta_{1}^{3}})\\
&=&a_{1}r^{1/3}s^{-2/3}(\sum_{m_{1},
n_{1}, \beta_{1}^{1},\beta_{1}^{2},\beta_{1}^{3}}b_{m_{1}, n_{1},
\beta_{1}^{1},\beta_{1}^{2},
\beta_{1}^{3}}K_{1}^{m_{1}}K_{2}^{n_{1}}E_{1}^{\beta_{1}^{1}}E_{2}^{\beta_{1}^{2}}E_{3}^{\beta_{1}^{3}})K_{1}.
\end{eqnarray*}

We also have the following
\begin{eqnarray*}
&&a_{2}K_{2}(\sum_{m_{1},n_{1},\beta_{1}^{1},\beta_{1}^{2},
\beta_{1}^{3}}b_{m_{1}, n_{1}, \beta_{1}^{1}, \beta_{1}^{2}, \beta_{2}^{3}}
K_{1}^{m_{1}}K_{2}^{n_{1}}E_{1}^{\beta_{1}^{1}}E_{2}^{\beta_{1}^{2}}E_{3}^{\beta_{1}^{3}})\\
&=&a_{2}r^{-1/3}s^{-1/3}(\sum_{m_{1},n_{1}, \beta_{1}^{1},\beta_{1}^{2},\beta_{1}^{3}}b_{m_{1}, n_{1},
\beta_{1}^{1},\beta_{1}^{2},
\beta_{1}^{3}}K_{1}^{m_{1}}K_{2}^{n_{1}}E_{1}^{\beta_{1}^{1}}E_{2}^{\beta_{1}^{2}}E_{3}^{\beta_{1}^{3}})K_{2}.
\end{eqnarray*}

Via detailed calculations and simplifications, we have the following
\[
\beta_{1}^{1}+\beta_{1}^{3}=1, \quad \beta_{1}^{2}+\beta_{1}^{3}=0.
\]
Based on the fact that $\beta_{i}^{j}$ are nonnegative integers, we have that 
\[
\beta_{1}^{1}=1,\quad \beta_{1}^{2}=0=\beta_{1}^{3}.
\]

Similarly, we can also have the following 
\[
\beta_{2}^{1}=\beta_{2}^{3}=0,\quad \beta_{2}^{2}=1.
\] 
Thus we have proved {\bf Case 1}.

{\bf Case 2}: Suppose that $\theta(K_{1})=a_{1}K_{2}$ and $\theta(K_{2})=a_{2}K_{1}$, we need to prove that
$\theta(E_{1})=b_{1}K_{1}^{m_{1}}K_{2}^{n_{1}}E_{2}$ and $\theta(E_{2})=b_{2}K_{1}^{m_{2}}K_{2}^{n_{2}}E_{1}$. 
Since the proof goes the same as in {\bf Case 1}, we will not repeat the detail here.
\qed 

Furthermore, we can easily verify that $E_{1}, E_{2}$ can not be exchanged. In particular, we  have the following result
\begin{cor}
Let $\theta \in Aut_{\C}(\V)$ be an algebra automorphism. Then for $l=1,2$, we have 
\[
\theta(K_{l})=a_{l}K_{l},\, \theta(E_{l})=b_{l}K_{1}^{m_{l}}K_{2}^{n_{l}}E_{l}
\]
where $a_{l}, b_{l} \in \C^{\ast}$ and $m_{l},n_{l} \in \Z$.
\end{cor}
\qed

Now we are going to prove one of the main results of this paper, which describes the algebra automorphism group of $\V$:
\begin{thm}
Let $\theta \in Aut_{\C}(\V)$ be an algebra automorphism of $\V$. Then for $l=1,2$, we have the following
\[
\theta(K_{l})=a_{l}K_{l},\quad \theta(E_{1})=b_{1}K_{1}^{a}K_{2}^{b}E_{1}, \quad \theta(E_{2})=b_{2}K_{1}^{c}K_{2}^{d}E_{2}
\]
where $a_{l}, b_{l} \in \C^{\ast}$ and $a,b,c,d\in \Z$ such that $b=c, a+b+d=0$.
\end{thm}
{\bf Proof:} Let $\theta$ be an algebra automorphism of $\V$ and suppose that 
\[
\theta(E_{1})=b_{1}K_{1}^{a}K_{2}^{b}E_{1},\quad
\theta(E_{2})=b_{2}K_{1}^{c}K_{2}^{d}E_{2}.
\]

Then we have the following
\begin{eqnarray*}
(K_{1}^{a}K_{2}^{b}E_{1})(K_{1}^{a}K_{2}^{b}E_{1})(K_{1}^{c}K_{2}^{d}E_{2})
&=& (r^{-1/3}s^{2/3})^{a} (r^{1/3}s^{1/3})^{b}(r^{-1/3}s^{2/3})^{2c}\\
& &(r^{1/3}s^{1/3})^{2d}K_{1}^{2a+c}K_{2}^{2b+d}E_{1}^{2}E_{2}\\
&=& r^{(-a+b-2c+2d)/3}s^{(2a+b+4c+2d)/3}\\
&& K_{1}^{2a+c}K_{2}^{2b+d}E_{1}^{2}E_{2};
\end{eqnarray*}

and
\begin{eqnarray*}
(K_{1}^{a}K_{2}^{b}E_{1})(K_{1}^{c}K_{2}^{d}E_{2})(K_{1}^{a}K_{2}^{b}E_{1})
&=& (r^{-1/3}s^{2/3})^{c} (r^{1/3}s^{1/3})^{d}(r^{-1/3}s^{-1/3})^{a}\\
&&(r^{-1/3}s^{2/3})^{a}(r^{-2/3}s^{1/3})^{b}(r^{1/3}s^{1/3})^{b}\\
&&K_{1}^{2a+c}K_{2}^{2b+d}E_{1}E_{2}E_{1}\\
&=& r^{(-2a-b-c+d)/3}s^{(a+2b+2c+d)/3}\\
&& K_{1}^{2a+c}K_{2}^{2b+d}E_{1}E_{2}E_{1};\\
\end{eqnarray*}

and
\begin{eqnarray*}
(K_{1}^{c}K_{2}^{d}E_{2})
(K_{1}^{a}K_{2}^{b}E_{1})(K_{1}^{a}K_{2}^{b}E_{1})&=&
(r^{-1/3}s^{-1/3})^{a} (r^{-2/3}s^{1/3})^{b}(r^{-2/3}s^{1/3})^{a}\\
&& (r^{-1/3}s^{2/3})^{b}K_{1}^{2a+c}K_{2}^{2b+d}E_{2}E_{1}^{2}\\
&=& r^{(-3a-3b)/3}s^{3b/3}K_{1}^{2a+c}\\
&&K_{2}^{2b+d}E_{2}E_{1}^{2}.
\end{eqnarray*}

Applying the automorphism $\theta$ to the first two-parameter quantum Serre relation
\[
E_{1}^{2}E_{2}-(r+s)E_{1}E_{2}E_{1}+rsE_{2}E_{1}^{2}=0,
\]
we have the following system of equations
\begin{eqnarray*}
-a+b-2c+2d&=& -2a-b-c+d;\\
-3a-3b&=&-2a-b-c+d;\\
2a+b+4c+2d&=& a+2b+2c+d;\\
3b&=&a+2b+2c+d.
\end{eqnarray*}

It is easy to see that the previous system of equations is reduced to the following system of equations
\begin{eqnarray*}
a+2b-c+d&=&0;\\
a-b+2c+d &=&0.
\end{eqnarray*}

Similarly, from the second two-parameter quantum Serre relation
\[
E_{1}E_{2}^{2}-(r+s)E_{2}E_{1}E_{2}+rsE_{2}^{2}E_{1}=0,
\] 
we also have the same system of equations as follows
\begin{eqnarray*}
a+2b-c+d&=&0;\\
a-b+2c+d &=&0.
\end{eqnarray*}

Solving the system
\begin{eqnarray*}
a+2b-c+d&=&0;\\
a-b+2c+d &=&0;
\end{eqnarray*}
we have that $b=c$ and $a+b+d=0$. Thus we have proved the theorem as desired.
\qed

\subsection{Hopf algebra automorphisms of $\V$}

In this subsection, we further determine all the Hopf algebra automorphisms of the Hopf algebra $\V$. Let us denote by $Aut_{Hopf}(\V)$ the group of all Hopf algebra automorphisms of $\V$. 

First of all, we have the following result 
\begin{thm}
Let $\theta \in Aut_{Hopf}(\V)$. Then for $l=1,2$, we have the following
\[
\theta(K_{l})=K_{l},\quad \theta(E_{l})=b_{l}E_{l},
\]
where $b_{l}\in \C^{\ast}$. In particular, we have 
\[
Aut_{Hopf}(\V)\cong (\C^{\ast})^{2}.
\]
\end{thm}
{\bf Proof:} Let $\theta \in Aut_{Hopf}(\V)$ be a Hopf algebra
automorphism of $\V$, then we have $\theta \in Aut_{\C}(\V)$. Therefore, 
we have the following
\begin{eqnarray*}
\theta(K_{l})=a_{l}K_{l};\\
\theta(E_{1})=b_{1}K_{1}^{a}K_{2}^{b}E_{1};\\
\theta(E_{2})=b_{2}K_{1}^{c}K_{2}^{d}E_{2};
\end{eqnarray*} 
where $a_{l}, b_{l}\in \C^{\ast}$ for $l=1,2$ and $a,b,c,d \in \Z$ such that $b=c$ and $a+b+d=0$. 

First of all, we need to prove that $a_{l}=1$ for $l=1,2$. Since $\theta$ is a
Hopf algebra automorphism, we have the following
\[
(\theta \otimes \theta)(\Delta(K_{l}))=\Delta(\theta(K_{l}))
\]
for $l=1,2$, which imply the following
\[
a_{l}^{2}=a_{l}
\]
for $l=1,2$. Thus we have $a_{l}=1$ for $l=1,2$ as desired. 

Second of all, we need to prove that $a=b=c=d=0$. Note that we have the following
\begin{eqnarray*}
\Delta(\theta(E_{1}))&=&\Delta(b_{1}K_{1}^{a}K_{2}^{b}E_{1})\\
&=&\Delta(b_{1}K_{1}^{a}K_{2}^{b})\Delta(E_{1})\\
&=&b_{1}(K_{1}^{a}K_{2}^{b}\otimes K_{1}^{a}K_{2}^{b}) (E_{1}\otimes 1+K_{1}^{2}K_{2}^{-1})\\
&=& b_{1}K_{1}^{a}K_{2}^{b}E_{1}\otimes K_{1}^{a}K_{2}^{b}+b_{1}K_{1}^{a}K_{2}^{b}K_{1}^{2}K_{2}^{-2}\otimes K_{1}^{a}K_{2}^{b}E_{1}\\
&=&\theta(E_{1})\otimes K_{1}^{a}K_{2}^{b}+K_{1}^{a}K_{2}^{b}K_{1}^{2}K_{2}^{-1}\otimes \theta(E_{1})
\end{eqnarray*}
and
 
\begin{eqnarray*}
(\theta\otimes \theta)(\Delta(E_{1}))&=&(\theta\otimes \theta)(E_{1}\otimes 1+K_{1}^{2}K_{2}^{-1}\otimes E_{1})\\
&=& \theta(E_{1})\otimes 1+\theta(K_{1}^{2}K_{2}^{-1})\otimes \theta(E_{1})\\
&=& \theta(E_{1})\otimes 1+ K_{1}^{2}K_{2}^{-1}\otimes \theta(E_{1}).
\end{eqnarray*} 
Since $\Delta(\theta(E_{1}))=(\theta\otimes \theta)\Delta(E_{1})$, we have $a=b=0$. Since $b=c$ and $a+b+d=0$, we have $a=b=c=d=0$.

In addition, it is obvious that the algebra automorphism $\theta$  defined by $\theta(K_{l})=K_{l}$ and $\theta(E_{l})=b_{l}E_{l}$ for $l=1, 2$ is a Hopf algebra automorphism of $\V$ as well. Therefore, we have proved the theorem.
\qed 

\section{Derivations and the first Hochschild cohomology group of $\U$}

In this section, we determine all the derivations of $\U$. In
particular, we prove each derivation of $\U$ can be uniquely written
as the sum of an inner derivation and a linear combination to
certain specifically defined derivations. As a result, we are able to 
prove that the first Hochschild cohomology group of $\U$ is 
a two-dimensional vector space over the base field $\C$. All 
these will be done through an embedding of the algebra $\U$ 
into a quantum torus, whose derivations had been described 
in \cite{OP}. This method has also been successfully used to 
compute the derivations of the algebra $U_{q}({\mathfrak sl_{4}^{+}})$ in 
\cite{LL}.

\subsection{The embedding of $\U$ into a quantum torus}

In this subsection, we embed the algebra $\U$ into a quantum torus,
which enables us to extend the derivations of $\U$ to derivations of 
the quantum torus. Note that the algebra $\U$ has a Goldie quotient 
ring, which we shall denote by $Q(\U)$. Inside the Goldie quotient 
ring $Q(\U)$ of $\U$, let us define the following new variables
\[
T_{1}=E_{1}, \quad T_{2}=E_{2}-\frac{1}{r-s}E_{3}E_{1}^{-1}, \quad
T_{3}=E_{3}.
\]

Concerning the relationships between the variables $T_{i}, i=1,2,3$, it
is easy to see that we have the following proposition
\begin{prop}
The following identities hold:
\begin{enumerate}
\item $T_{1}T_{2}=sT_{2}T_{1}$;\\
\item $T_{1}T_{3}=rT_{3}T_{1}$;\\
\item $T_{2}T_{3}=r^{-1}T_{3}T_{2}$.
\end{enumerate}
\end{prop}
\qed

Let us denote by $A^{3}$ the subalgebra of $Q(\U)$ generated by $T_{1}^{\pm
1}, T_{2}, T_{3}$, then we have the following
\begin{prop}
The subalgebra $A^{3}$ is the same as the subalgebra of $Q(\U)$ generated
by $E_{1}^{\pm 1}, E_{2}, E_{3}$. In particular, $A^{3}$ is a free
module over the subalgebra generated by $E_{2}, E_{3}$.
\end{prop}
\qed

Furthermore, let us denote by $A^{2}$ the subalgebra of $Q(\U)$ 
generated by $T_{1}^{\pm 1}, T_{2}, T_{3}^{\pm 1}$. Then we have 
the following proposition
\begin{prop}
The subalgebra $A^{2}$ is the same as the subalgebra of $Q(\U)$ generated
by $E_{1}^{\pm 1}, E_{2}, E^{\pm 1}_{3}$.
\end{prop}
\qed

Similarly, let us denote by $A^{1}$ the subalgebra of $Q(\U)$
generated by $T_{1}^{\pm 1}, T_{2}^{\pm 1}, T_{3}^{\pm 1}$. Thanks to 
{\bf Proposition 3.1}, we know that the indeterminates $T_{1}, T_{2}, T_{3}$ 
generate a quantum torus, which we shall denote by $Q_{3}=\C_{r,s}[T_{1}^
{\pm 1}, T_{2}^{\pm 1}, T_{3}^{\pm }]$. In particular, we have the
following
\begin{prop}
The algebra $A^{1}=Q_{3}=\C_{r,s}[T_{1}^{\pm 1}, T_{2}^{\pm 1},
T_{3}^{\pm 1}]$ is a quantum torus.
\end{prop}
\qed

Now let us define a linear map 
\[
\mathcal{I} \colon \U \longrightarrow Q_{3}
\]
from $\U$ into $A^{1}=Q_{3}$ as follows
\[
\mathcal{I}(E_{1})=T_{1}, \quad
\mathcal{I}(E_{2})=T_{2}+\frac{1}{r-s}T_{3}T_{1}^{-1},\quad
\mathcal{I}(E_{3})=T_{3}.
\]

It is easy to see that the linear map $\mathcal{I}$ can be extended to 
an algebra monomorphism from $\U$ into $A^{1}=Q_{3}$. Furthermore, it 
is straightforward to prove the following result:
\begin{thm}
Let us set $A^{4}=\U, \Sigma_{4}=\{T_{1}^{i} \mid i \in \Z_{\geq 0}\},
\Sigma_{3}=\{T_{3}^{i}\mid i \in \Z_{\geq 0}\},
\Sigma_{2}=\{T_{2}^{i}\mid i\in \Z_{\geq 0}\}$, then we have the following
\begin{enumerate}
\item $A^{3}=A^{4}\Sigma_{4}^{-1}$;\\
\item $A^{2}=A^{3}\Sigma_{3}^{-1}$;\\
\item $A^{1}=A^{2}\Sigma_{2}^{-1}$;\\
\item $A^{4}\subset A^{3}\subset A^{2}\subset A^{1}$;\\
\item The center of $A^{i}$ is the base field $\C$ for $i=1,2,3,4$.
\end{enumerate}
\end{thm}
\qed

From the reference \cite{OP}, one knows that a derivation $D$ of
the quantum torus $A^{1}=Q_{3}$ is of the form $D=ad_{t}+\delta$ where 
$ad_{t}$ is an inner derivation determined by some $t\in A^{4}$, and 
$\delta$ is a central derivation which acts on the variables 
$T_{i}, i=1, 2, 3$ as follows:
\[
\delta(T_{i})=\alpha_{i} T_{i}
\]
for $\alpha_{i}\in \C$.

Let $D$ be a derivation of $\U=A^{4}$. According to the previous
theorem, one can extend the derivation $D$ to a derivation of $A^{i}$
for $i=3, 2, 1$, and thus a derivation of the quantum torus
$A^{1}$. We still denote the extension by $D$. Therefore, as 
a derivation of the quantum torus $A^{1}$, the derivation 
$D$ can be decomposed as follows
\[
D=ad_{t}+\delta
\]
where $ad_{t}$ is an inner derivation determined by some $t\in A^{1}$,
and $\delta$ is a central derivation of $A^{1}$
which is defined by $\delta(T_{i})=\alpha_{i}T_{i}$ for $\alpha_{i}
\in \C, i=1,2,3$. 

Now we are going to prove that the element $t$ can actually be chosen
from the algebra $\U$ and the scalars $\alpha_{1},\alpha_{2}$ and 
$\alpha_{3}$ are related to each other. In particular, we have the 
following key lemma.
\begin{lem}
The following is true:
\begin{enumerate}
\item The element $t$ can be chosen from $\U$;\\
\item We have $\alpha_{3}=\alpha_{1}+\alpha_{2}$;\\
\item We have $D(E_{i})=ad_{t}(E_{i})+\alpha_{i}E_{i}$ for $i=1,2,3$.
\end{enumerate}
\end{lem} 
{\bf Proof:} First of all, we show that $t\in A^{2}$. Suppose that we 
have $t=\sum_{i,j, k}a_{i,j,k}T_{1}^{i}T_{3}^{j}T_{2}^{k}\in A_{1}$. 

If $k\geq 0$ for all $k$, then we are done. Otherwise, let us set
\[
t_{-}=\sum_{k<0}a_{i,j,k}T_{1}^{i}T_{3}^{j}T_{2}^{k}
\]
and 
\[
t_{+}=\sum_{k\geq 0}a_{i,j,k}T_{1}^{i}T_{3}^{j}T_{2}^{k}.
\]
Then we have that $t=t_{-}+t_{+}$. 

First of all, we have the following
\begin{eqnarray*}
D(T_{1}) &=& ad_{t}(T_{1})+\delta(T_{1})\\
&=& (t_{-}T_{1}-T_{1}t_{-})+(t_{+}T_{1}-T_{1}t_{+})+\alpha_{1}(T_{1})
\end{eqnarray*}
for some $\alpha_{1} \in \C$.

Since $D$ is a derivation of the algebra $\U$ and the variable $T_{1}$ 
is in the algebra $A^{4}=\U$, we have that the element $D(T_{1})$ is in 
the algebra $A^{4}$, and furthermore in the algebra $A^{2}$. Note that
the elements of $A^{2}$ don't involve negative powers of the variable 
$T_{2}$, thus we have the following
\[
t_{-}T_{1}=T_{1}t_{-}. 
\]

Therefore, we have the following
\begin{eqnarray*}
T_{1}(\sum_{k<0}a_{i,j,k}T_{1}^{i}T_{3}^{j}T_{2}^{k})&=& 
(\sum_{k<0}a_{i,j,k}r^{j}s^{k}T_{1}^{i}T_{3}^{j}T_{2}^{k})T_{1}\\
&=& (\sum_{k<0}a_{i,j,k}T_{1}^{i}T_{3}^{j}T_{2}^{k})T_{1}.
\end{eqnarray*}

Thus we have $k=0$. Since we are supposed to have $k<0$ in the
expression of $t_{-}$, we have run into a contradiction. Therefore,
we have that $t_{-}=0$ and $t=t_{+}\in A^{2}$. Similarly, we can also 
prove that $t\in A^{3}$.

Since the algebra $A^{3}$ is also generated by the elements $E_{1}^{\pm 1}, E_{2},
E_{3}$, we have the following 
\[
t=\sum_{i,j\geq 0, k\geq
0}a_{i,j,k}E_{1}^{i}E_{3}^{j}E_{2}^{k}.
\]

Applying the derivation $D$ to the variable $E_{3}$, we can further
prove that $i\geq 0$. Therefore, we have proved that $t\in A^{4}=\U$ 
as desired.

Since we have $D=ad_{t}+\delta$ for some $t\in \U$, we have 
\begin{eqnarray*}
D(E_{2})&=&ad_{t}(E_{2})+\delta(E_{2})\\
&=& (tE_{2}-E_{2}t)+\delta(T_{2}+\frac{1}{r-s}T_{3}T_{1}^{-1})\\
&=& (tE_{2}-E_{2}t)+\alpha_{2}(T_{2}+\frac{1}{r-s}T_{3}T_{1}^{-1})\\
&&+\frac{1}{r-s}(\alpha_{3}-\alpha_{1}-\alpha_{2})T_{3}T_{1}^{-1}\\
&=& (tE_{2}-E_{2}t)+\alpha_{2}E_{2}+\frac{1}{r-s}(\alpha_{3}-\alpha_{1}-\alpha_{2})T_{3}T_{1}^{-1}.
\end{eqnarray*}

Note that $D(E_{2})\in \U$, we have
$\frac{1}{r-s}(\alpha_{3}-\alpha_{1}-\alpha_{2})T_{3}T_{1}^{-1}\in
A^{4}$. Thus we have the following
\[
\alpha_{3}=\alpha_{1}+\alpha_{2}
\]

and
\[
D(E_{2}))=ad_{t}(E_{2})+\alpha_{2}E_{2}.
\]

So we have the proved the lemma as desired.
\qed

Now let us define two derivations $D_{1}, D_{2}$ of the algebra $\U$ 
as follows:
\begin{eqnarray*}
D_{1}(E_{1})=E_{1},\quad D_{2}(E_{2})=0,\quad D_{1}(E_{3})=E_{3};\\
D_{2}(E_{1})=0,\quad D_{2}(E_{2})=E_{2},\quad D_{2}(E_{3})=E_{3}.
\end{eqnarray*}

Based on the previous lemma, we have the following
\begin{thm}
Let $D$ be a derivation of $\U$. Then we have 
\[
D=ad_{t}+\mu_{1}D_{1}+\mu_{2}D_{2}
\]
for some $t\in \U$ and $\mu_{i}\in \C$ for $i=1,2$.
\end{thm}
\qed

Recall that the Hochschild cohomology group in degree 1 of $\U$ is
denoted by $HH^{1}(\U)$, which is defined as follows
\[
HH^{1}(\U)\colon = Der(\U)/InnDer(\U).
\]
where $InnDer(\U)\colon= \{ad_{t} \mid t \in \U \}$ is the Lie algebra of 
inner derivations of $\U$. It is well known that $HH^{1}(\U)$ is a 
module over $HH^{0}(\U) \colon= Z(\U)=\C$. 

Now we state the structural result for the first Hochschild cohomology
of $\U$.

\begin{thm}
The following is true:
\begin{enumerate}
\item Every derivation $D$ of $\U$ can be uniquely written as follows:
\[
D=ad_{t}+\mu_{1}D_{1}+\mu_{2}D_{2}
\]
where $ad_{t}\in InnDer(\U)$ and $\mu_{1}, \mu_{2}\in \C$.\\
\item The first Hochschild cohomology group $HH^{1}(\U)$ of $\U$ is a 
two-dimensional vector space spanned by $\overline{D_{1}}$ 
and $\overline{D_{2}}$.
\end{enumerate}
\end{thm}
{\bf Proof:} Suppose that we have
$ad_{t}+\mu_{1}D_{1}+\mu_{2}D_{2}=0$, then we need to prove that 
$\mu_{1}=\mu_{2}=ad_{t}=0$. Let us set
$\delta=\mu_{1}D_{1}+\mu_{2}D_{2}$. Then $\delta$ is a derivation of the algebra $\U$. 

Note that we can extend the derivation $\delta$ to a derivation of $A^{1}$, and we also
have $ad_{t}+\delta=0$ as a derivation of $A^{1}$. Furthermore, we
have the following
\[
\delta(T_{1})=\mu_{1}T_{1},\quad \delta(T_{2})=\mu_{2}T_{2}, \quad \delta(T_{3})=(\mu_{1}+\mu_{2})T_{3}.
\]

Thus the derivation $\delta$ is a central derivation of the quantum
torus $A^{1}$. According to the result in \cite{OP}, we have that 
$ad_{t}=0=\delta$. Thus we have $\mu_{1}=\mu_{2}=0$ as desired. 

So we have proved the uniqueness for the decomposition of $D$, which 
further implies the second statement of the theorem.
\qed

\end{document}